\documentclass[12pt,a4paper]{article}

\usepackage[T1]{fontenc}
\usepackage[utf8]{inputenc}
\usepackage{authblk}
\usepackage{amsfonts}
\usepackage{mathtools}
\usepackage{slashed}
\usepackage{amsmath,amssymb}
\usepackage{graphicx,xcolor}
\usepackage{cite}
\usepackage{hyperref}

\usepackage{alphabeta} 

\definecolor{refcol}{rgb}{0.9,0.1,0.1}
\hypersetup{colorlinks=true,linkcolor=blue,citecolor=refcol,urlcolor=cyan,linktocpage}
  
\usepackage[capitalize]{cleveref}
\usepackage{tikz}
\usetikzlibrary{shapes}
\usetikzlibrary{plotmarks}

\usepackage{setspace}

\usepackage[left=2.5cm,right=2.5cm,top=2.5cm,bottom=2.5cm]{geometry}
\allowdisplaybreaks
\newcommand{\be}{\begin{equation}}
\newcommand{\ee}{\end{equation}}
\newcommand{\bea}{\begin{eqnarray}}
\newcommand{\eea}{\end{eqnarray}}

\def\XXint#1#2#3{{\setbox0=\hbox{$#1{#2#3}{\int}$ }
\vcenter{\hbox{$#2#3$ }}\kern-.6\wd0}}

\newcommand{\Qp}{\mathbb{Q}_p}
\newcommand{\Zp}{\mathbb{Z}_p}
\newcommand{\bupsi}{\mbox{\bf \textpsi}}  

\begin{document}
\begin{titlepage}

\title{
{\huge\bf A $p$-arton Model for Modular Cusp Forms}\\ 
}

\bigskip\bigskip\bigskip\bigskip\bigskip

\author{{\bf Parikshit Dutta}${}^{1}$\thanks{{\tt parikshitdutta@yahoo.co.in}}\hspace{8pt} and \hspace{2pt}  
                    {\bf Debashis Ghoshal}${}^2$\thanks{{\tt d.ghoshal@gmail.com, dghoshal@mail.jnu.ac.in}} \\  
\hfill\\
${}^1${\it Asutosh College, 92 Shyama Prasad Mukherjee Road,}\\
{\it Kolkata 700026, India}
\hfill\\              
${}^2${\it School of Physical Sciences, Jawaharlal Nehru University,}\\
{\it New Delhi 110067, India}
}

\date{%
%
\bigskip\bigskip
\begin{quote}
\centerline{{\bf Abstract}}
{\small 
We propose to associate to a modular form (an infinite number of) complex valued functions on the $p$-adic numbers $\Qp$ for each
prime $p$. We elaborate on the correspondence and study its consequences in terms of the Mellin transform and the $L$-function
related to the form. Further we discuss the case of products of Dirichlet $L$-functions and their Mellin duals, which are convolution 
products of $\vartheta$-series. The latter are intriguingly similar to non-holomorphic Maass forms of weight zero as suggested by their
Fourier coefficients.  
}
\end{quote}
}

\bigskip

\end{titlepage}
\thispagestyle{empty}\maketitle\vfill \eject

\tableofcontents

\section{Introduction}\label{sec:Intro} 
The modular groups SL(2,$\mathbb{Z}$) and its subgoups, while of interest to mathematicians for their myriad manifestations, also play 
an important role in models of statistical mechanics where they are related to the partition function viewed as a function of complexified 
temperature and other external parameters. Holomorphic modular forms, with definite transformation properties under a subgroup of the
modular group, are essential ingredients in conformally invariant quantum field theories in two dimensions, and more specifically in string 
theory and black hole physics, where both holomorphic and non-holomorphic modular forms arise (see e.g., 
\cite{Dabholkar:2019wzh,Dabholkar:2012zz} and references therein). 
The usual holomorphic modular forms are defined as $q$-series expansion in which the coefficients satisfy certain multiplicative properties. 
The coefficients of a modular cusp form define an associated $L$-series, the two being related by the Mellin transform. Thanks to the
multiplicative properties of its coefficients, an $L$-series admits an Euler product over prime numbers.  

Another class of $L$-series is the family of Dirichlet $L$-series, of which the Riemann zeta function is the most well known member 
\cite{hmedwards}, associated to the multiplicative Dirichlet characters. Each of these also admits a Euler product representation 
in terms of prime numbers. The product form allows one to connect the local factors to the $p$-adic number fields $\Qp$. We have 
recently constructed \cite{Dutta:2020qed} pseudodifferential operators (as generalisations of the notion of the Vladimirov derivative in 
\cite{VVZ1994p}) which incorporate the Dirichlet characters and shown that the $L$-series can be expressed as the trace of an operator
in an appropriate vector space. These operators act on the Bruhat-Schwarz class of locally constant complex valued square integrable 
functions on $\Qp$. The traces of these operators are well defined, being restricted to operators that act on the subspace $L^2(p^{-1}\Zp)$ 
of functions with compact support in $p^{-1}\Zp\subset\Qp$. In a larger collaboration \cite{Chattopadhyay:2018bzs} we have also initiated 
a programme to construct a random ensemble of unitary matrix models (UMM) for the prime factors of the Riemann zeta function, and then 
combine them to get a UMM for the latter. We find that the partition function of the UMM can be written as a trace of the Vladimirov derivative 
restricted to $L^2(p^{-1}\Zp)$. In this approach the Riemann zeta function is essentially a partition function in the sense of statistical 
mechanics---see also \cite{Spector1990,Julia1990,Bakas1991,Julia1994,DG:FZRZ} for `physical models' of $L$-functions, a key difference is 
that we use the orthonormal basis provided by the wavelets \cite{Kozyrev:2001} in $\Qp$. In \cite{Dutta:2020qed} we also constructed 
pseudodifferential operators corresponding to $L$-series of modular cusp forms. We showed that a family of locally constant functions, 
called the Kozyrev wavelets \cite{Kozyrev:2001}, known to be the eigenfunctions of the Vladimirov derivative, is the set of common 
eigenfunctions of all these pseudodifferential operators.
 
The purpose of this article is to propose an association between a modular cusp form and complex valued functions in $L^2(\Qp)$, one function
for each prime $p$. More precisely, the correspondence is between the Fourier expansion of the cusp forms and the functions on $\Qp$. In other 
words, we argue that a cusp form is equivalent to a vector in $\otimes_p L^2(\Qp)$. This decomposition reminds us of the \emph{parton} model 
of hadrons (which, like primes, begins fortuitously with p). A $p$-adic Mellin transform of these vectors, when combined for all primes, is shown 
to be related to the $L$-series corresponding to the cusp form. We also discuss Hecke operators in terms of the raising-lowering operators on 
the wavelet basis. This requires the definition of an appropriate inner product at the level of the $p$-artons. We examine two possibilities and 
study their properties. In an attempt to elaborate on these, we propose to define a class of toy $L$-functions (by taking a product of two Dirichlet 
$L$-functions) that mimic the Euler product form of $L$-functions associated with modular forms. As a pleasant surprise, we find that the objects 
associated to these bear an intriguing resemblance to non-analytic Maass forms \cite{MaassZL}. It is likely that these are indeed toy examples 
of Maass-like forms, however, we have only been able to verify their behaviour under the $\mathrm{Im }(z)\to -1/\mathrm{Im }(z)$ transform of 
the modular group.      
 
In the following, we review a few relevant facts concerning holomorphic cusp forms of the modular groups related to SL(2,$\mathbb{Z}$) 
(\cref{sec:MFLfn}) and wavelets on the $p$-adic numbers (\cref{sec:BaseQp}) that also help us establish the notation used. In 
\cref{sec:Parton} we propose to associate complex valued functions on $\Qp$, one for each prime $p$, which we call the $p$-artons, 
to a cusp form, and discuss their Mellin transforms. Realisation of the Hecke operators in this description is taken up in \cref{sec:Hecke}, 
where we propose possible inner products on the spaces of $p$-artons and their Mellin duals. Finally, in \cref{sec:ProdDLfn} we study 
the modular objects associated to products of two Dirichlet $L$-functions and point out their relation to non-analytic Maass forms. We 
conclude in \cref{sec:Summary} with a brief summary and a comment on the holographic nature of the correspondence proposed in this 
paper. 

\section{Modular forms and associated $L$-functions}\label{sec:MFLfn}
The discrete subgroup\footnote{More precisely, the relevant groups are the projective special linear groups PSL(2, $\mathbb{R}$) =
SL(2,$\mathbb{R}/\{\pm\}$) and PSL(2,$\mathbb{Z}$) = SL(2,$\mathbb{Z})/\{\pm 1\}$.} 
$\mathrm{SL}(2,\mathbb{Z}) = 
\left\{ \gamma = \left(\begin{array}{cc}
a & b\\
c & d\end{array}\right)\;\big|\: a, b, c, d \in\mathbb{Z},\, ad-bc=1\right\}$ of the special linear group $\mathrm{SL}(2,\mathbb{R})$, is the symmetry 
group of lattices $\Lambda$ in the complex plane $\mathbb{C}$ \cite{Serre:Course,Koblitz:ECMF,Warner:Crash,MIT:MFLF}. It is sometimes called 
the full modular group, and has the following \emph{congruence subgroups} of finite indices
\begin{align}
\Gamma_0(N) &= \left\{\gamma\in \mathrm{SL}(2,\mathbb{Z})\, |\, c \equiv 0\: \text{mod}\: N\right\}\nonumber\\  
\Gamma_1(N) &= \left\{\gamma\in \mathrm{SL}(2,\mathbb{Z})\, |\, a,d \equiv 1\:\text{and}\: c \equiv 0\: \text{mod}\: N\right\}
\label{CongruentSubGp}\\
\Gamma(N) &= \left\{\gamma\in \mathrm{SL}(2,\mathbb{Z})\, |\, a,d \equiv 1\:\text{and}\: b,c \equiv 0\: \text{mod}\: N\right\}\nonumber
\end{align}     
(Since the conditions are empty for $N=1$, $\Gamma(1)$ is the full modular group $\mathrm{SL}(2,\mathbb{Z})$.) They are ordered
as follows: $\Gamma(N) \subset \Gamma_1(N)\subset \Gamma_0(N) \subset \Gamma(1)$. Of these, $\Gamma(N)$, called the 
principal congruence subgroup of level $N$, is the kernel of the homomorphism $\mathrm{SL}(2,\mathbb{Z}) \rightarrow 
\mathrm{SL}(2,\mathbb{Z}/N\mathbb{Z})$. 

If $\Gamma$ is a discrete subgroup such that $\Gamma(N) \subset \Gamma \subset\Gamma(1)$, where $N$ is the smallest such integer, it 
is referred to as a congruence subgroup of level $N$. The natural action of the group $\mathrm{SL}(2,\mathbb{R})$ on the upper half plane 
$\mathbb{H} = \{z \,:\, \mathrm{Im }(z) > 0\}$ restricts to that of $\Gamma$, which partitions it into equivalence classes. A \emph{fundamental 
domain} $\mathcal{F}$ is a subset of $\mathbb{H}$ representing a $\Gamma$-equivalence class. For example, the fundamental domain of the 
full modular group $\Gamma(1)$ is $\left\{z\in\mathbb{H} \,|\, -\frac{1}{2} \le \mathrm{Re }(z) \le \frac{1}{2}, |z| \ge 1 \right\}$ (ignoring some
double counting of points on the boundaries). 

A \emph{modular form} \cite{Koblitz:ECMF} $f : \mathbb{H}\to\mathbb{C}$ of \emph{weight} $k$ and \emph{level} $N$ associated to a Dirichlet 
character\footnote{A Dirichlet character (modulo $N$) is a group homomorphism $\chi_N\in \mathrm{Hom}(G(N),\mathbb{C}^*)$ from the 
multiplicative group $G(N)=\left(\mathbb{Z}/N\mathbb{Z}\right)^*$ of invertible elements of $\mathbb{Z}/N\mathbb{Z}$ to $\mathbb{C}^*$. It is a 
multiplicative character. It customarily extended to all integers by setting $\chi_N(m)=0$ for all $m$ which share common factors with $N$ 
\cite{Serre:Course}.} $\chi_N$ modulo $N$, is a holomorphic form on the upper half plane $\mathbb{H}$ that transforms, under the action of a 
discrete subgroup $\Gamma(N) \subset \Gamma \subset \Gamma(1)$, as
\begin{equation}
f(\gamma z) \equiv f\left(\frac{az + b}{cz + d}\right) = \chi_N(d) (cz + d)^{k} f(z), \quad \label{ModTrans}
\end{equation}
where $k\in\mathbb{N}$, $N\in\mathbb{Z}$, $\gamma \in \Gamma\subset \mathrm{SL}(2,\mathbb{Z})$. The modular form vanishes identically 
unless $k$ is an even integer. 

Using $z\to z+1$ in \cref{ModTrans}, one sees that a modular form of the full modular group $\Gamma(1)$ (i.e., of level 1) is a
periodic function in $z$, hence it has the following Fourier expansion ($q$-expansion) in $q=e^{2\pi i z}$
\begin{equation}
f(z) = \sum_{n=0}^\infty a(n) q^n = \sum_{n=0}^\infty a(n) e^{2\pi i n z} \label{qExpand} 
\end{equation}
A \emph{cusp form} is a modular form that vanishes as $\mathrm{Im }(z)\to i\infty$, or equivalently at $q=0$. Thus $a(0)=0$ for a cusp form. 
It is conventional, and often convenient, to normalise the first coefficient $a(1)$ of a cusp form to 1, which is what we shall assume in what 
follows. An equivalent description of modular forms is in terms of scaling functions on $\mathbb{C}/\Lambda$, where $\Lambda$ is a lattice 
left invariant by the action of a subgroup of the modular group. 

Modular forms of weight $k$ form a \emph{finite dimensional} complex vector space $M_k\left(\Gamma(1)\right)$ and of these the subset of 
cusp forms is a subspace $S_k\left(\Gamma(1)\right)$. Similar notions exist for modular forms of congruence subgroups $\Gamma \subset 
\Gamma(1)$ of level $N$. However, a cusp form of a congruence subgroup $\Gamma$ is required to vanish as $z$ approaches certain rational 
points on $\mathbb{R} = \partial\mathbb{H}$ (equivalently, at certain point on the unit circle $|q|=1$), in addition to $z\to i\infty$ ($q=0$). These 
additional points in the fundamental domain are images of  $\mathrm{Im }(z) \to i\infty$. 

The Dirichlet series of a cusp form $f=\sum_n a(n) q^n$ is defined by the coefficients in its $q$-expansion:  
\begin{equation}
L(s,f) =  \sum_{n=1}^\infty \frac{a(n)}{n^s} = 1 + \frac{a(2)}{2^{s}} + \frac{a(3)}{3^{s}} + \frac{a(4)}{4^{s}} + \cdots
\label{ModularL}
\end{equation}
where we have used the normalisation $a(1)=1$. This series converges uniformly to a holomorphic function of $s$ in a half-plane to the right 
of $\mathrm{Re }(s) = \sigma+1$ as long as the coefficients $|a(n)|$ are bounded by some power $n^\sigma$. The corresponding $L$-function 
associated to the cusp form $f$ is then defined by an analytic continuation to the complex $s$-plane. For a cusp form of weight $k$, the series 
above converges in $\mathrm{Re }(s) > \frac{k+1}{2}$. The series is also related to the cusp form as
\begin{equation}
L(s,f) = \mathcal{M}[f](s) = \frac{(2\pi)^s}{\Gamma(s)} \int_0^\infty dy\, y^{s-1} f(iy)  \label{FMellinL}
\end{equation}
i.e., it is the Mellin transform of $f(iy)$.

The discriminant function $\Delta(z) = \left(2\pi\right)^{12} \displaystyle{\prod_{n=1}^\infty \left(1-q^n\right)}$ is an example of a cusp
form of weight 12 (and level 1) of the full modular group. Ramanujan noticed that the coefficients in its $q$-expansion 
\begin{equation}
\Delta(z)  =  \sum_{n=1}^\infty \tau(n) q^n
\end{equation}
satisfy the following properties
\begin{align}
\begin{split}
\tau(mn) &= \tau(m) \tau(n) \; \text{if gcd }(m,n)=1\\
\tau(p^{m+1}) &= \tau(p) \tau(p^m) - p^{11} \tau(p^{m-1}) \;\, \text{for } p \:\; \text{a prime and } m>0 \\
|\tau(p)| &\le 2 p^{11/2}
\end{split}\label{RamanujanTau}
\end{align}
The function $\tau : \mathbb{N}\rightarrow \mathbb{Z}$ is known as the Ramanujan $\tau$-function. (The first two properties were proved 
by Mordell, while the proof for the bound was provided by Deligne.) 
More generally, the coefficients $a(n)$ of a modular form of weight $k$ and level $N$ satisfy
\begin{align}
\begin{split}
a(mn) &= a(m) a(n) \; \text{if gcd }(m,n)=1\\
a(p^{m+1}) &= a(p)\, a(p^m) - \chi(p) p^{k-1} a(p^{m-1}) \;\, \text{for } p \:\; \text{a prime and } m>0 
\end{split}\label{Recursion}
\end{align}
The coefficient function $a: \mathbb{N}\rightarrow \mathbb{C}$ is said to define a multiplicative character \cite{Serre:Course}. The convergence 
of the series also puts a bound on the growth of the coefficients, which for a cusp form is $|a(p)| \le Cp^{{(k-1)}/{2}}$ for $C$ of order 1. Due to 
the above properties of its coefficients, the $L$-function of a cusp form $f$ \cref{ModularL} too admit an Euler product form
\begin{equation}
L(s,f) = \sum_{n=1}^\infty \frac{a(n)}{n^s}  = \prod_{p\,\in\,{\mathrm{primes}}} \frac{1}{\left(1 - a(p) p^{-s} +\chi(p) 
p^{k-1} p^{-2s}\right)} \label{EulerModL}   
\end{equation}  
This is analogous to the Dirichlet $L$-functions \cref{DirichL}, however, unlike that case, each local factor $L_p(s,f) = \left(1 - a(p) p^{-s} 
+ \chi(p) p^{k-1} p^{-2s}\right)^{-1} $ in the denominator of $L(s,f)$ is a quadratic in $p^{-s}$. 

\section{A wavelet basis for complex valued functions on $\Qp$}\label{sec:BaseQp}
The prime factors in the Euler product of Dirichlet $L$-functions \cref{DirichL}, as well as $L$-functions associated with cusp forms can be 
related to complex valued functions on the $p$-adic space $\Qp$. Using the ultrametric $p$-adic norm, which itself is an example of a 
complex valued function $|\cdot|_p:\Qp\to\mathbb{C}$, one may write the prime factor $(1-p^{-s})^{-1}$ in the Riemann zeta function as
\begin{equation}
\frac{1}{1-p^{-s}} \equiv \zeta_p(s) = \frac{p}{p-1} \int_{\Zp} dx\, |x|_p^{s-1},\qquad\mathrm{Re }(s) > 0 
\label{LocZeta}
\end{equation}
where the $dx$ is the (complex valued) translation invariant Haar measure on $\Qp$. The integral may also be thought
of as
\begin{equation}
\frac{p}{p-1} \int_{\Zp^\times} d^\times x\, |x|_p^s,\qquad\mathrm{Re }(s) > 0
\label{LocZMellin} 
\end{equation}
in which the integral is over the non-zero elements in $\Zp$ with respect to the multiplicative (scale invariant) measure 
$d^\times x = dx/|x|_p$ of $\Qp^\times$ viewed as a multiplicative group. The integrand $|x|_p^s$, a multiplicative character 
on $\Qp$, is another complex valued function.

Let us introduce two other functions
\begin{align}
\begin{split}
\chi_p(x) &= e^{2\pi i x},\quad x\in\Qp\\
\Omega_p(x,x_0) & = \begin{cases}
1 & \mathrm{if } |x-x_0|_p \le 1\\
0 & \text{otherwise}
\end{cases},\quad
x,x_0\in\Qp
\end{split}
\label{AddCharIndFn}
\end{align}
where the first one is an additive character and the second is an indicator function for a unit ball centred at $x_0$
\cite{Gelfand1968representation,VVZ1994p}. This additive character defines the Fourier transform $\mathcal{F}[f]$ of a function 
$f:\Qp\to\mathbb{C}$ as $\mathcal{F}[f](k) = \int_{\Qp} dx\, \chi_p(kx) f(x)$ and its inverse $\mathcal{F}^{-1}$. The indicator function 
$\Omega_p(x)$ retains its form after Fourier transform, i.e., $\mathcal{F}[\Omega_p](k) = \Omega_p(k)$. In this sense, it is an analogue 
of the Gaussian function $e^{-\pi x^2}$ on $\mathbb{R}$.  These are the ingredients of a family of orthonormal functions\cite{Kozyrev:2001} 
on $\Qp$
\begin{align}
\psi^{(p)}_{n,m,j}(x) &= p^{-\frac{n}{2}} e^{\frac{2\pi i}{p} j p^{n} x} \Omega_p(p^n x - m) \\
\int_{\Qp} dx\, \psi^{(p)}_{n,m,j}(x) \psi^{(p)}_{n',m',j'}(x) & = \delta_{nn'} \delta_{mm'} \delta_{jj'}
\label{KozyWvlet}
\end{align}
where $n\in\mathbb{Z}$, $m\in\Qp/\Zp$ and $j=1,2,\cdots, p-1$. These functions, which we shall refer to as the {\em Kozyrev wavelets}, provide 
an orthonormal basis for the set $L^2(\Qp)$ of square-integrable Bruhat-Schwarz (locally constant) functions on $\Qp$. These wavelets may be
considered to be analogous to the (generalised) Haar wavelets on $\mathbb{R}$. All these functions can be obtained from the {\em mother 
wavelets} $\psi_{0,0,j}(x)$ by the action of the affine group `$ax+b$' with suitable choices for $a$ and $b$.

The usual definition of the derivative of a function does not work due to the totally disconnected topology of $\Qp$. However, a
pseudodifferential operator, called the generalised Vladimirov derivative \cite{VVZ1994p,WilsonBook}, is defined by an integral 
kernel
\begin{align}\label{VladD}
\begin{split}
D_{(p)}^\alpha f(x) &= \frac{1}{\Gamma_{(p)}(-\alpha)}\, \int dx'\, \frac{f(x') - f(x)}{|x' - x|_p^{\alpha+1}}\\
\text{where,}\quad\Gamma_{(p)}(-\alpha) &= \int_{\Qp^\times} \frac{dx}{|x|_p} \, e^{2\pi ix} |x|_p^{-\alpha} 
\end{split}
\end{align}
for $\alpha$ in a suitable half plane of $\mathbb{C}$ and elsewhere by analytic continuation when appropriate. The Kozyrev 
wavelets are eigenfunctions of the Vladimirov derivatives 
\begin{equation}
D_{(p)}^\alpha \psi^{(p)}_{n,m,j} (x) = p^{\alpha(1-n)} \psi^{(p)}_{n,m,j} (x) \label{VDonKozy}
\end{equation}
corresponding to the eigenvalues $p^{\alpha(1-n)}$. The operator $\log_p D_{(p)}$ can be defined in the limit $\alpha\to 0$ as 
\begin{equation*}
\log_p D_{(p)} = \displaystyle{\lim_{\alpha\to 0}}\frac{D_{(p)}^\alpha - 1}{\alpha\ln p}
\end{equation*}
such that $\log_p D_{(p)} \psi^{(p)}_{n,m,j} (x) = (1-n) \psi^{(p)}_{n,m,j} (x)$. These wavelets are also common eigenfunctions of a more general 
class of Vladimirov derivatives twisted by a multiplicative character\cite{Dutta:2020qed}.

Since our interest will be in the eigenvalues, which depend only on the quantum number $n$ related to scaling (and not $m$ and 
$j$ related to translation and phase), we shall restrict our attention to the set of eigenfunctions
\begin{equation}
\psi_{n,0,1}^{(p)} (x) = p^{-\frac{n}{2}} \chi_p(p^{n-1} x)\, \Omega_p( p^n x)\quad \longleftrightarrow \quad\left| 1 - n
\right\rangle_{(p)}  \label{waveket}
\end{equation}
where we have used the alternative ket-vector notation labelled by the eigenvalue of the wavelet.  

Let us also define the raising/lowering operators $a^{(p)}_\pm$ on the wavelets
\begin{equation}
a^{(p)}_\pm\psi^{(p)}_{n,0,1}(x)  = \psi^{(p)}_{n\pm 1,0,1}(x) \quad\longleftrightarrow\quad  a^{(p)}_\pm | n \rangle_{(p)} 
= | n \mp 1\rangle_{(p)} \label{raise:lower}
\end{equation}
that changes the scaling quantum number by one. With these operators, together with $\log_p D_{(p)}$, one can define 
$J^{(p)}_\pm = a^{(p)}_\pm\log_p D_{(p)}$ and $J^{(p)}_3 = \log_p D_{(p)}$ the commutator algebra of which generate an 
$\mathfrak{sl}(2,\mathbb{R})$ symmetry of the wavelets \cite{Dutta:2018qrv}. 

We shall be more interested in the subset spanned by the wavelets
\begin{equation}
\left\{ |m\rangle_{(p)} \; | \; m=0, 1, 2, \cdots\right\} \;\buildrel{n=1-m}\over{\longleftrightarrow}\; \left\{ \psi^{(p)}_{n,0,1}(x)\; | \; 
n=1,0,-1,-2,\cdots\right\} \label{eq:subset}
\end{equation}  
to define the subspace $\mathcal{H}^{(p)}_- \subset L^2(p^{-1}\Zp)$. It consists of wavelets supported on the compact subset $p^{-1}\Zp$ 
and the corresponding eigenvalues of the operator $D_{(p)}$ on these basis functions are $\left\{1, p, p^2, \cdots\right\}$, i.e., positive integer 
powers of $p$. When restricted to this subspace, we shall demand that the lowering operator $a^{(p)}_+$ annihilates the mother wavelet 
corresponding to the lowest eigenvalue of the Vladimirov derivative\footnote{By an abuse of notation, we shall continue to use the same symbol 
for the restrictions of the raising-lowering to $\mathcal{H}^{(p)}_-$. However, since it is the subspace that is of primary interest to us, this should 
hopefully not be a cause of confusion.}
\begin{equation}
a^{(p)}_+ | 0\rangle_{(p)} = 0\; \longleftrightarrow\; a^{(p)}_+ \psi^{(p)}_{1,0,1}(x) = 0\label{annihilate}
\end{equation}
Thus the wavelet $\psi^{(p)}_{1,0,1}(x)$ is like a `ground state' or `lowest weight state' in this subspace, and the other wavelets 
arise from repeated applications of the raising operator $a^{(p)}_-$ on it, namely $|n\rangle =\big(a^{(p)}_-\big)^n |0\rangle_{(p)}$.

\section{Vectors for modular (cusp) forms}\label{sec:Parton}
We are now ready to propose an association between a modular cusp form and functions on $\otimes_p\Qp$, more precisely on
$\otimes\mathcal{H}^{(p)}_-$ , i.e., those spanned by the wavelets in the subspace \cref{eq:subset}. For definiteness, we shall discuss the 
association for cusp forms of the principal congruent subgroup $\Gamma(N)$ of the modular group. 

First we use prime factorisation of a natural number $n\in\mathbb{N}$ to relate it to a wavelet in $\otimes\mathcal{H}^{(p)}_-\subset
\otimes_p\Qp$ as follows 
\begin{align}
n &= \prod_{p\in\text{Primes}} p^{n_p}\quad\longmapsto \nonumber\\ 
\bigotimes_{p} \left| n_p\right\rangle_{(p)} &= |n_2\rangle_{(2)} \otimes |n_3
\rangle_{(3)} \otimes |n_5\rangle_{(5)} \otimes  |n_7\rangle_{(7)} \otimes \cdots\label{PrimeFact}
\end{align}
Since all but a finite number of $n_p$s are zero, the associated wavelets are the lowest weight or ground states $|0\rangle_{(p)}$ 
corresponding to the wavelets $\psi^{(p)}_{0,0,1}(x)$ for the prime $p$. Only a finite number of wavelets are, therefore, {\em 
non-trivial} in this correspondence. Similar association between natural numbers as vectors in the Hilbert space of fictitious 
quantum systems (dubbed arithmetic gas or primon gas) using the prime factorisation has been made in 
Refs.\cite{Spector1990,Julia1990,Bakas1991,Julia1994}. The present proposal relates to a mathematically well-defined space of
functions, however, without any apparent physical meaning. On the other hand, the types of modular forms and $L$-functions for
which we shall apply this correspondence are much larger.  

Recall the expansion of a cusp form in \cref{qExpand,ModularL} in which the leading coefficient $a(1)$ has been normalised to 
unity. To this cusp form, we associate the following vector in $\otimes_p\mathcal{H}^{(p)}_-$
\begin{align}
f = \sum_{n=1}^\infty a(n) q^n &= \sum_{n_p=0}^\infty \left(\prod_{p} a(p^{n_p})\right)  q^{\prod p^{n_p}}\nonumber\\
\longmapsto\quad |\mathfrak{f}\rangle &= \sum_{n_p=0}^\infty \bigotimes_{p} a(p^{n_p}) \left| n_p\right\rangle_{(p)}\nonumber\\
&=  \sum_{{n_2,n_3,n_5,\atop\cdots = 0}}^\infty\!\!\! a(2^{n_2}) |n_2\rangle_{(2)} \otimes a(3^{n_3}) |n_3\rangle_{(3)} 
\otimes a(5^{n_5}) |n_5\rangle_{(5)} \otimes \cdots\nonumber\\
&= \sum_{n_p=0}^\infty \bigotimes_{p} a(p^{n_p}) a_-^{n_p} |0\rangle_{(p)} \label{qExpPFactorised}
\end{align}
where we have used the multiplicative property of the coefficients \cref{Recursion} for arguments which are coprime to each other. 
It is understood that for $n_p=0$, the coefficients $a(1)=1$ for all $p$. We note that in the above too, most of the terms are `trivial', 
i.e., most vectors correspond to the \emph{ground state} $|0\rangle_{(p)} = |n_p=0\rangle \in \mathcal{H}^{(p)}_-$. In the last line, we 
use the fact that ket vectors of higher \emph{occupation number} $n_p \ge 1$ can be obtained from the ground state by the action of 
the raising operator \cref{raise:lower}, i.e., $|n_p\rangle = a_-^{n_p} |0\rangle_{(p)}$. 

The vector $|\mathfrak{f}\rangle\in\otimes_p\mathcal{H}^{(p)}_-$ appears to be a complicated one which is really \emph{entangled}.
However, thanks to the multiplicative property of the coefficients \cref{Recursion}, it simplifies to a product form
\begin{equation}
|\mathfrak{f}\rangle = \sum_{n_p=0}^\infty \bigotimes_{p} a(p^{n_p}) a_-^{n_p} |0\rangle_{(p)} 
= \left(1 - a(p) a_- + p^{k-1}\chi(p) a_-^2\right)^{-1} |0\rangle_{(p)}\equiv \bigotimes_{p}\, \left|\mathfrak{f}_{(p)}\right\rangle 
\label{ModularAtom}
\end{equation}
as can be seen by expanding the right hand side and comparing it with the left hand side. We would like to think of the vector
$\left|\mathfrak{f}_{(p)}\right\rangle\in\mathcal{H}^{(p)}_-$  as the $p$-th $p$-arton, i.e., the `part' of the cusp form $f$ at the 
prime $p$. It is interesting to note that the operator that acts on the `ground state' to generate the $p$-arton resembles the form of 
the local $L$-function at a prime $p$ in \cref{EulerModL}. Indeed, that was a clue to find the factorisation above.

All this can equivalently be expressed in terms of the complex valued wavelet functions, if we introduce the \emph{coordinate basis}, 
familiar in quantum mechanics, consisting of generalised kets $\left\{ \left|x_{(p)}\right\rangle , x_{(p)}\in\Qp \right\}$ which satisfy the 
\emph{orthogonormality condition} 
\begin{equation*}
\left\langle x_{(p)} \big| x'_{(p')}\right\rangle = \delta_{pp'}\, \delta\!\left(x_{(p)} - x'_{(p)}\right)
\end{equation*}
involving the Dirac $\delta$-functions. In this basis the components $\langle x_{(p)}|\mathfrak{f}_{(p)}\rangle$ (also called the 
{\em wavefunction}) of $|\mathfrak{f}_{(p)}\rangle$ are
\begin{equation}
f_{(p)}(x_{(p)}) \equiv \left\langle x_{(p)} \big| \mathfrak{f}_{(p)}\right\rangle = \sum_{n_p=0}^\infty a(p^{n_p}) \left\langle x_{(p)} 
\big| n_p \right\rangle
= \sum_{n_p=0}^\infty a(p^{n_p}) \psi^{(p)}_{1-n_p,0,1}(x_{(p)}) \label{ModFormPart}
\end{equation}
i.e., a linear combination of wavelets compactly supported on $p^{-1}\Zp\subset\Qp$. 

\begin{quote}
\emph
{\noindent
In summary, according to this correspondence, a cusp form $f : \mathbb{H} \to \mathbb{C}$ is equivalent to the infinite set of functions 
$f_{(p)} : \Qp\to\mathbb{C}$, one function for each prime $p$. The two are equivalent in the sense that from $f$ we can get 
$\left(f_{(2)}, f_{(3)}, f_{(5)}, \cdots\right)$ and vice versa.}
\end{quote}     

We would like to show that the association proposed in \cref{qExpPFactorised,ModularAtom,ModFormPart} has other ramifications. 
To that end, we shall consider a Mellin transformation of the wavelets. There are different proposals for $p$-adic Mellin transforms 
\cite{GoldfeldHundley,Gelfand1968representation}, however, all involve integrands with multiplicative characters on $\Qp$. The
following definition\footnote{In \cite{Gelfand1968representation} the multiplicative character in the kernel is taken to be unitary, which is 
satisfied if $s$ in \cref{pMellin} is purely imaginary. The definition in \cite{GoldfeldHundley} (see definitions 2.8.4 and 2.8.5) uses 
a \emph{normalised unitary character} $\omega$, which has a \emph{conductor} $p^N$, which for this case is $N=1$.} that we shall 
employ is similar to the one in Ref.\cite{GoldfeldHundley} 
\begin{equation}
\tilde{g}_\omega(s) \equiv \mathcal{M}_{(p,\omega)}[g](s) = \int_{\Qp^\times} d^\times\! x \: e^{\frac{2\pi i\ell}{p}x |x|_p}\, |x|_p^{s}\, g(x),
\qquad s\in\mathbb{C} \text{ and } \ell=0,1,\cdots,p-1
\label{pMellin}
\end{equation}
where $d^\times\! x = dx/|x|_p$ is the scale invariant measure and the kernel contains the unitary character $\omega_\ell(x) = 
e^{\frac{2\pi i\ell}{p} x |x|_p}$. The character, which determines a phase depending on the leading $p$-adic `digit' of $x$, has 
been mentioned in \cite{VVZ1994p}. The transformation satisfies the scaling property
\begin{equation*}
\mathcal{M}_{(p,\omega)}\left[g(\alpha x)\right](s) = |\alpha|_p^{-s} \mathcal{M}_{(p,\omega)}[g(x)](s)   
\end{equation*}
similar to usual Mellin transforms, with $\omega_\ell \rightarrow \omega_{\ell\alpha|\alpha|_p}$, which is also a primitive $p$-th 
root of unity, if $\omega_\ell$ is. In particular, for $\alpha=p^n$ the unitary character does not change, consequently 
$\mathcal{M}_{(p,\omega)}\left[g(p^n x)\right](s) = p^{ns} \mathcal{M}_{(p,\omega)}[g(x)](s)$ may be used to relate the Mellin 
transform of the wavelets. The inverse Mellin transform
\begin{equation}
\mathcal{M}^{-1}_{(p,\omega)}[\tilde{g}_\omega](x) = \sum_{\ell=0}^{p-1} e^{-\frac{2\pi i\ell}{p} x |x|_p}\, \frac{\ln p}{2\pi} 
\int_0^{2\pi/\ln p} \!\!dt\, |x|_p^{-it}\, \tilde{g}_\omega(it) \label{InvpMellin}
\end{equation}
involves a discrete Fourier transform of the character $\omega$ and the path of integration is from $0$ to $2\pi/\ln p$ along $t$, the imaginary
axis in the complex $s$-plane.

With this definition, the Mellin transform of the Kozyrev wavelet $\psi_{n,0,1}(x)$ is
\begin{equation}
\mathcal{M}_{(p,\omega)}[\psi_{n,1,0}](s) = -\left(\frac{1}{p(1 - p^{-s})}-\frac{1}{p^s - 1}\delta_{\ell,0} - \delta_{\ell,p-1} \right) 
p^{n\left(s-\frac{1}{2}\right)} \label{KozyMellin}
\end{equation}
One can verify that the inverse transform of the above is the Kozyrev wavelet by an explicit calculation. It is relevant to point out that even 
though the contribution of the last term to the inverse transform adds to zero, while the contributions from the first two terms reproduce the 
wavelet function, the discrete Fourier transform involving the character $\omega$ is absolutely crucial for the inverse transform to work. As 
an aside, it is interesting to note that the Mellin transform of a Kozyrev wavelet is related to the eigenvalue of the generalized Vladimirov 
derivative $D^{-s}$.

{}From \cref{ModFormPart} and \cref{KozyMellin}, we get 
\begin{equation}
\mathcal{M}_{(p,\omega)}\left[\mathfrak{f}_{(p)}(x_{(p)})\right](s) =  c_p(\ell,s) \sum_{n_p=0}^\infty a(p^{n_p})\,
p^{(1-n_p)\left(s-\frac{1}{2}\right)} \label{MellinFAtom}
\end{equation}
where 
\begin{equation}
c_p(\ell,s) = - \frac{1}{p(1 - p^{-s})} + \frac{\delta_{\ell,0}}{p^s - 1} + \delta_{\ell,p-1} = \begin{cases}
-p^{-s}\Gamma_p(s), &\ell = 0\\
p^{-1}\zeta_p(s), &\ell = 1,2,\cdots,p-2\\
p^{-1}\zeta_p(s) - 1, &\ell = p-1  
\end{cases}
\label{MTCoeff}
\end{equation}
is the $n_p$-independent factor in the parenthesis in \cref{KozyMellin}. Now we can combine the results for all primes: the Mellin transform 
of the `wavefunction' of \cref{ModFormPart} is
\begin{align}
\mathcal{M}_{(p,\omega)} \left[ \big\langle (\xi_{(2)},\xi_{(3)},\xi_{(5)},\cdots) \big| \mathfrak{f}\big\rangle \right] (s) 
&= \prod_p \mathcal{M}_{(p,\omega)}\left[\mathfrak{f}_{(p)}(\xi_{(p)})\right] (s) \nonumber\\
{} &= \prod_p c_p(\ell,s) \sum_{n_p=0}^\infty a(p^{n_p}) p^{(1-n_p)\left(s-\frac{1}{2}\right)} \nonumber\\ 
{} &= \left(\prod_p c_p(\ell,s) p^{s-\frac{1}{2}}\right)\, L\left(s-\frac{1}{2},f\right) \label{AdMellin}
\end{align}
Curiously the $L$-function we get this way has its argument shifted from $s$ to $s-\frac{1}{2}$. The infinite product in the prefactor also 
depends on $\ell$, which in turn depends on $p$.

\section{Raising, lowering and Hecke operators}\label{sec:Hecke}
It would be instructive to understand the proposed decomposition of a modular form in terms of $p$-adic wavelets in the context of the 
profound Hecke theory of modular forms. However, this is outside the scope of our present understanding. We shall need to know how to 
extend the definition of the Kozyrev wavelets, which are defined on $\Qp$, to the projective space $\mathbb{P}(\Qp) \sim \Qp\cup\{\infty\}$ 
so that their transformation under the full group GL(2,$\Qp$) can be addressed. 
For now, we shall set a more modest goal of an operational understanding of the Hecke operators $T(m)$, $m\in\mathbb{N}$ and their 
algebraic properties. 

Recall that the Hecke operators $T(m)$, $m\in\mathbb{N}$, are a set of commuting operators whose action on the modular
form  $f(z)=\sum a(n) e^{2\pi i nz}$ is to return the coefficients in the $q$-expansion as eigenvalues
\cite{Serre:Course,Koblitz:ECMF,GoldfeldHundley,MIT:MFLF,Warner:Crash}
\begin{equation}
T(m) f(z) = a(m) f(z) \label{HeckeOnModForm}
\end{equation}
In other words a modular form is an eigenvector of the Hecke operators with the eigenvalues as the coefficients in its $q$-expansion. 
They satisfy
\begin{align}
\begin{split}
T(m) T(n) &= T(mn) \text{ for } m\nmid n\\
T(p) T(p^\ell) &= T(p^{\ell+1}) + \chi(p) p^{k-1} T(p^{\ell-1}) 
\end{split}
\label{HeckeAlg}
\end{align}
Alternatively, the action of the Hecke operator $T(n)$ on the underlying lattice can be understood as a sum of sublattices of index $n$: 
$T(n) \Lambda =\! \displaystyle{\sum_{[\Lambda:\Lambda']=n}}\! \Lambda'$. 

The Hecke operator $T(m)$ can be written as a sum of two operators\cite{Koblitz:ECMF,Serre:Course}, the first of which, $V(m)$ gives a 
new series by replacing each  $q^n$ in $f$ by $q^{mn}$ 
\begin{equation} 
\left(V(m)f\right) (z) = \sum_{n=1}^\infty a(n) q^{mn} = \sum_{n=1}^\infty a(n) e^{2\pi i mnz} = f(mz) \label{HeckeVOp}
\end{equation}
while the action of the second, $U(m)$, is to define a new series by keeping only those terms $q^n$ that are divisible by $m$
\begin{equation} 
\left(U(m)f\right) (z) = \sum_{{n=1 \atop (m\mid n)}}^\infty a(n) q^{\frac{n}{m}}  
= \frac{1}{m}\, \sum_{j=0}^{m-1} f\left(\frac{z+j}{m}\right) \label{HeckeUOp} 
\end{equation}
Note that $U(m) V(m) = \mathbf{1}$, but since $V(m) U(m)$ deletes terms in \cref{qExpand} not divisible by $m$, $V(m) U(m) \ne 
\mathbf{1}$. A Hecke operator for prime argument can be written as
\begin{equation}
T(p) = U(p) + \chi(p) p^{k-1} V(p) \label{HeckePrime}
\end{equation}
Thus, if we denote $\left(T(p)f\right) (z) = \sum b(n) q^n$, then $b(n) = a(pn) + \chi(p) p^{k-1} a(n/p)$.

The actions of the $U$ and $V$ operators remind us of the raising and lowering operators \cref{raise:lower} on the wavelets. Consequently
\begin{align*}
\sum_{n=0}^\infty a(p^n) |n\rangle\: & \buildrel a_{-}\over\longrightarrow\: \sum_{n=0}^\infty a(p^n) |n+1\rangle 
= \sum_{n=1}^\infty a(p^{n-1}) |n\rangle\\    
\sum_{n=0}^\infty a(p^n) |n\rangle\: & \buildrel a_{+}\over\longrightarrow\: \sum_{n=1}^\infty a(p^n) |n-1\rangle 
= \sum_{n=0}^\infty a(p^{n+1}) |n\rangle\\    
{} & \qquad\qquad = a(p) \sum_{n=0}^\infty a(p^n) |n\rangle - \chi(p) p^{k-1} \sum_{n=1}^\infty a(p^{n-1}) |n\rangle  
\end{align*}
Hence, we propose to define the following Hecke operator for a prime argument
\begin{equation}
T(p) = a_{+} + \chi(p) p^{k-1} a_{-} \label{pHecke}
\end{equation}
to obtain 
\begin{equation*}
T(p) \left| \mathfrak{f}_{(p)}\right\rangle = a(p) \left| \mathfrak{f}_{(p)}\right\rangle
\end{equation*}
It is easy to see that this will lead to the correct eigenvalues \cref{HeckeOnModForm}.

\subsection{Inner product I}\label{ssec:IPone}
The vector space of modular (cusp) forms of weight $k$ is equipped with the Petersson inner product, defined as  
\begin{equation}
\langle f | g \rangle = \int_{\mathcal{F}} d^2z\, \left(\mathrm{Im}(z)\right)^{k-2} f^*(z) g(z)
\label{PeterssonIP}
\end{equation}
which is invariant under the action of the modular group. The Hecke operators, in general, are not hermitian, rather they satisfy
\begin{equation*}
T^\dagger(p) = \chi^*(p) T(p) \label{HeckeConjugation}
\end{equation*}
where, hermitian conjugation is defined in terms of inner product \cref{PeterssonIP}. In other words, $\chi^{*\frac{1}{2}} T$ is a hermitian 
operator. 

How does $T(p)$ in \cref{pHecke} acting on the space of vectors $\left\{\left| \mathfrak{f}_{(p)}\right\rangle\right\}$ behave under conjugation?
First, one needs to define an inner product. At this point we do not have an understanding of the transformation properties of the wavelets 
under the {\em local modular group} GL(2,$\Qp$), therefore, we do not know how to define a GL(2,$\Qp$) invariant inner product analogous to
the one in \cref{PeterssonIP}. We shall examine two choices, both seemingly natural, in turn. The first is defined in analogy with the Petersson 
inner product above, however, using only the subgroup of scaling transformations, as follows:
\begin{equation}
\big\langle \mathfrak{f}_{(p)} \big| \mathfrak{g}_{(p)} \big\rangle = \int_{\Qp^\times} d^\times x \, \left| x\right|^{k-1} 
\mathfrak{f}_{(p)}^*(x)\, \mathfrak{g}_{(p)}(x) \label{pPeterssonIP}
\end{equation}
Recall that the effects of the operators $a_{\pm}$ on the wavelets is to scale $x$ by a factor of $p^{\pm 1}$ (as well as a change in the overall
normalisation such the raised/lowered wavelet functions remain orthonormal). By a change of variable of integration in the inner product defined 
above, it is easy to check, that
\begin{equation*}
a_{+}^\dagger = p^{k-1} a_{-} \quad\text{ and }\quad a_{-}^\dagger = p^{1-k} a_{+} 
\end{equation*}
Hence, it follows that $T^\dagger(p) = \chi^*(p) T(p)$, as well as that
\begin{equation*}
\left(a_{\mathfrak{g}}(p) - \chi^*(p) a_{\mathfrak{f}(p)}\right)\,\big\langle \mathfrak{f}_{(p)} \big| \mathfrak{g}_{(p)} \big\rangle = 0
\end{equation*}
Therefore, either $\mathfrak{f}_{(p)} = \mathfrak{g}_{(p)}$, in which case the expression in the parenthesis vanishes establishing
the hermiticity of $T(p)$, or the two vectors associated to distinct cusp forms are orthogonal $\langle \mathfrak{f}_{(p)} | 
\mathfrak{g}_{(p)} \rangle = 0$. 

We should hasten to add that the arguments (and hence the conclusions) above are formal. This is because the manipulation by scaling 
and redefinition of the variable of integration work only for functions in $L^2(\Qp)$. We do, however, need to restrict to a subspace 
$L^2(p^{-1}\Zp)$, in which scaling is not a symmetry as it may take us out of the subspace. Consequently, the relevant operators have
to be projected back to the subspace of our interest. It is not obvious that the process of conjugation will commute with the projection.
The arguments in \cref{ssec:IPtwo} based on another proposal for an inner product does not rely on scaling.

We would like to present an interesting observation, the analysis of the implication of which we shall leave for future. In the
inner product \cref{pPeterssonIP}, let us write the function $\mathfrak{g}_{(p)}(x)$ as the inverse (generalised) Mellin transform 
of its (generalised) Mellin transform, as defined in \cite{GoldfeldHundley}, to get
\begin{align}
\big\langle \mathfrak{f}_{(p)} \big| \mathfrak{g}_{(p)} \big\rangle &= \int_{\Qp^\times}\!\! d^\times x \, \left|x\right|^{(k-1)/2} 
\mathfrak{f}_{(p)}^*(x)\; \mathcal{M}_\omega^{-1}\left[\mathcal{M}_\omega\left[\left|x\right|^{(k-1)/2} \mathfrak{g}_{(p)}
\right]\right]\!(x) \nonumber\\
{} &= \int_{\Qp^\times}\!\! d^\times x \, \left|x\right|^{\frac{k-1}{2}} \mathfrak{f}_{(p)}^*(x)\!\! \sum_{{\omega\atop (\text{mod } p^N)}} 
\!\!\!\frac{\ln p}{2\pi} \int_0^{\frac{2\pi}{\ln p}}\!\!\! dt\, |x|_p^{-it} \omega^*(x) \mathcal{M}_\omega\left[\left| x\right|^{\frac{k-1}{2}} 
\mathfrak{g}_{(p)}\right](it) \nonumber\\
{} &= \sum_{{\omega\atop (\text{mod } p^N)}}\!\!\! \frac{\ln p}{2\pi} \int_0^{\frac{2\pi}{\ln p}}\!\!\!  dt \left(\!\mathcal{M}_\omega \left[ 
\mathfrak{f}_{(p)}\right]\!\left(\frac{k-1}{2} + it\right)\!\right)^* \mathcal{M}_\omega \left[ \mathfrak{g}_{(p)}\right]\!\left(\frac{k-1}{2} 
+ it\right)   \label{ParsevalType}
\end{align}
namely, a Parseval type 
It is interesting to note that the argument of the Mellin transform in the last line is $\frac{k}{2} - \frac{1}{2} + it$, which is shifted by half 
from the position of the conjectured zeroes of the $L$-function \cref{ModularL} on the line $\text{Re }(s) = \frac{k}{2} + it$. However, 
the poles of the local $L$-function in \cref{EulerModL} are likely to be on the line $\text{Re }(s) = \frac{k-1}{2} + it$. This can be
verified for those cases where the character $\chi$ is trivial, for example, the $L(s,\Delta)$ the coefficients of which are the 
Ramanujan $\tau$-function.    

\subsection{Inner product II}\label{ssec:IPtwo}
In order to discuss the other inner product and its applications, it turns out to be more convenient to modify the Kozyrev wavelets 
\cref{KozyWvlet} to define a set that is orthonormal with respect to the scale invariant measure $d^\times x$ on $\Qp^\times$. This
modification is described in the \cref{app:AppA}. We would also like to redefine the $p$-artonic modular forms \cref{ModFormPart}
by rescaling the coefficients with a weight dependent factor as follows  
\begin{equation}
\mathbf{f}_{(p)}(x_{(p)}) = \sum_{n_p=0}^\infty p^{-\frac{k-1}{2}n_p} a(p^{n_p}) \bupsi^{(p)}_{1-n_p,0,1}(x_{(p)}) 
\label{modModFormPart}
\end{equation}
The rescaling is motivated by the bound on the growth of the coefficients of cusp forms. The inner product we define is the simple 
overlap integral 
\begin{equation}
\big( \mathbf{f}_{(p)} \big| \mathbf{g}_{(p)} \big) = \int_{\Qp^\times} d^\times x \; \mathbf{f}_{(p)}^*(x)\, 
\mathbf{g}_{(p)}(x) \label{pSimpleIP}
\end{equation}
of the modified $p$-artonic wavefunctions.

We use the orthonormality conditions \cref{modOrthonorm} to express the action of the raising and lowering operators 
\cref{raise:lower,annihilate} on the function \cref{modModFormPart} to write 
\begin{align}
\mathbf{a}_{+} \mathbf{f}_{(p)} (x) &= \sum_{n_p=1}^\infty \int_{\Qp^\times} d^\times y\, \bupsi^{(p)}_{2-n_p,0,1}(x)
\bupsi^{(p)*}_{1-n_p,0,1}(y) \mathbf{f}_{(p)} (y) \nonumber\\
\mathbf{a}_{-} \mathbf{f}_{(p)} (x) &= \sum_{n_p=0}^\infty \int_{\Qp^\times} d^\times y\, \bupsi^{(p)}_{-n_p,0,1}(x)
\bupsi^{(p)*}_{1-n_p,0,1}(y) \mathbf{f}_{(p)} (y) 
\label{RLmodMod} 
\end{align}
{}From this one may check that $\mathbf{a}_{\pm}^\dagger = \mathbf{a}_{\mp}$ with respect to the inner product \cref{pSimpleIP} above. 
Moreover, using the same set of equations, we get
\begin{align*}
\mathbf{a}_+\mathbf{f}_{(p)}(x_{(p)}) &= p^{-\frac{k-1}{2}} a(p) \bupsi^{(p)}_{1,0,1}(x_{(p)}) + p^{-\frac{k-1}{2}} 
\sum_{n_p=1}^\infty p^{-\frac{k-1}{2}n_p} a(p^{n_p+1}) \bupsi^{(p)}_{1-n_p,0,1}(x_{(p)}) \\
\mathbf{a}_-\mathbf{f}_{(p)}(x_{(p)}) &= p^{\frac{k-1}{2}} \sum_{n_p=1}^\infty p^{-\frac{k-1}{2}n_p} a(p^{n_p-1}) 
\bupsi^{(p)}_{1-n_p,0,1}(x_{(p)})
\end{align*} 
Hence the following combination
\begin{equation}
\mathbf{T}(p)\mathbf{f}_{(p)} \equiv \left(\mathbf{a}^{(p)}_+ + \chi(p) \mathbf{a}^{(p)}_-\right) \mathbf{f}_{(p)} = p^{-\frac{k-1}{2}} 
a(p) \mathbf{f}_{(p)} \quad\text{ and }\quad \mathbf{T}^\dagger(p) = \chi^*(p) \mathbf{T}(p) \label{modHecke}
\end{equation}
behaves like the Hecke operator $T(p)$ ({\em cf.} \cref{HeckeOnModForm,HeckeConjugation}). Thanks to these equations
\begin{equation*}
\int_{\Qp^\times} d^\times x\, \mathbf{g}^*_{(p)} (x)\left(\mathbf{T}(p) - \chi^*(p)\mathbf{T}^\dagger(p)\right) \mathbf{f}_{(p)} (x) 
= 0
\end{equation*}
implies that
\begin{equation}
p^{-\frac{k-1}{2}}\left(a_{\mathbf{f}}(p) - \chi(p) a^*_{\mathbf{g}}(p)\right) \int_{\Qp^\times} d^\times x\, \mathbf{g}^*_{(p)} (x)
\mathbf{f}_{(p)} (x) = 0   \label{modOrtho}
\end{equation}
from which we conclude that $\mathbf{f}(p)$ and $\mathbf{g}(p)$ are orthogonal (with respect to the inner product \cref{pSimpleIP})
and that $\chi^{*\frac{1}{2}}(p) a_{\mathbf{f}}(p)$ are real.

Following the steps leading to the equality in \cref{ParsevalType} we find the corresponding Parseval type identity for this inner product
\begin{align}
\big( \mathbf{f}_{(p)} \big| \mathbf{g}_{(p)} \big) &= \int_{\Qp^\times}\!\! d^\times x\, \mathbf{f}_{(p)}^*(x)\; 
\mathcal{M}_\omega^{-1}\left[\mathcal{M}_\omega\left[ \mathbf{g}_{(p)}\right]\right]\!(x) \nonumber\\
{} &= \sum_{{\omega\atop (\text{mod } p^N)}}\!\!\! \frac{\ln p}{2\pi} \int_0^{\frac{2\pi}{\ln p}}\!\!\!  dt \left( \mathcal{M}_\omega \left[ 
\mathbf{f}_{(p)}\right]\!(it)\right)^* \mathcal{M}_\omega \left[ \mathbf{g}_{(p)}\right]\!(it)   \label{modParsevalType}\\
{} &= \frac{\ln p}{2\pi} \int_0^{\frac{2\pi}{\ln p}}\!\!\!  dt\, \left( L_{\mathbf{f}(p)}\Big(\frac{k-1}{2} + it\Big)\right)^* 
\left( L_{\mathbf{g}(p)}\Big(\frac{k-1}{2} + it\Big)\right) \nonumber
\end{align}
In arriving at the last step above we have used the following. In the Mellin transform $\mathcal{M}_\omega \left[\mathbf{f}_{(p)}\right]\!(it) = 
\mathbf{c}_p(\ell,it) p^{it} L_{\mathbf{f}(p)}\Big(\frac{k-1}{2} + it\Big)$, only the prefactor depends on $\omega(\ell)$, therefore, the sum
in the discrete Fourier transform yields 1 for the argument $s=it$ (see \cref{app:AppA}) and that the $t$-integral leads to a Kronecker
delta with which one of the sums is evaluated trivially. Thus, both the left and right hand sides reduce to $\sum_{n_p} p^{-(k-1)n_p} 
a^*_{\mathbf{f}}(p^n) a_{\mathbf{g}}(p^n)$. The left hand side vanishes if $\mathbf{f}$ and $\mathbf{g}$ arise from distinct modular forms. 
Then the above is an orthogonality condition for the corresponding modular $L$-functions. 

It should be instructive to check this by a direct computation, which we have unfortunately not been able to verify. Instead, we offer an indirect 
argument. Let us parametrise the `roots' of the denominator of the local function $L_p(s,f)$ \cref{EulerModL}, which is a quadratic in $p^{-s}$, 
as $a_1(p) = p^{(k-1)/2} e^{i\alpha_1(p)}$ and $a_2(p) = p^{(k-1)/2} e^{-i\alpha_2(p)}$ so that 
\begin{align*}
a(p) &= p^{\frac{k-1}{2}} \left(e^{i\alpha_1(p)} + e^{i\alpha_2(p)}\right) = 2\cos \frac{1}{2}(\alpha_1 + \alpha_2)\, p^{\frac{k-1}{2}} 
e^{\frac{i}{2} \left(\alpha_1 - \alpha_2\right)} \\
\chi(p) &= e^{i\left(\alpha_1(p) - \alpha_2(p)\right)}
\end{align*}
Notice that this is consistent with the condition on the growth of the coefficients $a(p)$. The function $L_p(s,f)$ is the  generating function 
\begin{equation}
\frac{1}{1 - 2t\cos\theta + t^2} = \sum_{n=0}^\infty U_n(\cos\theta) t^n \label{Cheby2GenFn}
\end{equation}
of the family of orthogonal Chebyshev polynomials of type II, denoted by $U_n(x)$, with $\theta = \frac{1}{2}(\alpha_1 + \alpha_2)$ and 
$t = p^{\frac{k-1}{2}-s} e^{\frac{i}{2} \left(\alpha_1 - \alpha_2\right)}$. The Chebyshev polynomials of type II satisfy the three-term recursion 
relation
\begin{equation}
U_{n+1}(\xi) = 2\xi\, U_n(\xi) - U_{n-1}(\xi),\qquad U_0(\xi) = 1\text{ and } U_1(\xi)=2\xi\label{Cheby2Recursion}
\end{equation}
In terms of trigonometric functions $U_n(\cos\theta) = {\sin((n+1)\theta)}/{\sin\theta}$. Furthermore, thanks to these recursion relation, they 
satisfy the multiplicative property \cref{Recursion}. The orthogonality of the local functions $L_p(s,f)$, and hence of the product functions 
$L(s,f)$ will follow as a consequence of the orthogonality of the Chebyshev polynomials. The appearance of the Chebyshev polynomials of 
type II in this context has been noticed\footnote{We became aware of these results after submitting a version of this article to the arXiv.} in 
Refs.~\cite{ConreyEtAl,SerreJAMS}. In the next section, we study a simpler class of functions, for which these properties will be manifest by 
construction.

\section{Products of Dirichlet $L$-functions and Maass-like forms}\label{sec:ProdDLfn}
Since the modular forms and the associated $L$-functions are rather complicated objects, we shall instead study a family of functions related 
to the Dirichlet $L$-function (corresponding to the Dirichlet character $\nu$)
\begin{equation}
L(s,\nu) = \sum_{n=1}^\infty \frac{\nu(n)}{n^s} = \prod_{p\,\in\,{\mathrm{primes}}} \frac{1}{\left(1 - \nu(p) p^{-s}\right)} 
\label{DirichL}
\end{equation}
which, we shall see, is simpler. We cannot use these $L$-functions directly because, unlike the $L$-functions associated with a modular cusp 
form \cref{EulerModL}, where the local factors at prime $p$, $\left(1 - a(p) p^{-s} + \chi(p) p^{k-1} p^{-2s}\right)^{-1}$ are quadratic functions of 
$p^{-s}$, the local factors $L_p(s,\nu) =\left(1 - \nu(p) p^{-s}\right)^{-1}$ above are linear. Therefore, in order to mimic the properties of a modular 
$L$-function, we consider a product of two Dirichlet $L$-functions to define the function ${}_2\mathsf{L}(s,\nu)$ as follows
\begin{align}
{}_2\mathsf{L}(s,\nu) =  L(s,\nu) L(s,\nu^*) &= \prod_p \frac{1}{\left(1 - \nu(p) p^{-s}\right)\left(1 - \nu^*(p) p^{-s}\right)} \nonumber\\
\sum_{n=1}^\infty \frac{a(n)}{n^s} &= \prod_p \frac{1}{\left(1 - 2\cos (\mathrm{arg}\,\nu_p) p^{-s} + p^{-2s}\right)}\label{augDirichL}
\end{align}
where $\nu^*$ is the complex conjugate of the character $\nu$ (it is further assumed that $\nu$ is not a principal character) and 
$a(n) = \displaystyle{\sum_{d|n}}\nu(d)\,\nu^*\!\left(\frac{n}{d}\right)$ is the convolution of the characters\cite{ApostolANT}. Notice that the function 
${}_2\mathsf{L}(s,\nu)$ is still meromorphic as a function of $s$. Formally it has the same form as \cref{EulerModL} with $a(p) = \nu(p)+\nu^*(p) 
= 2\cos (\mathrm{arg}\,\nu_p)$, $k=1$, and $\chi=1$ (the trivial character). 

The local factor ${}_2\mathsf{L}_p(s,\nu)$ at a prime $p$ can be recognised as the generating function \cref{Cheby2GenFn}, hence
\begin{equation}
{}_2\mathsf{L}_p(s,\nu) = \frac{1}{\left(1 - 2\cos (\mathrm{arg}\,\nu_p) p^{-s} + p^{-2s}\right)} = \sum_{n_p=0}^\infty U_{n_p}(\cos\xi) p^{-sn_p}
\label{GenCheby2}
\end{equation}
where $U_{n_p}(\xi)$ are the  Chebyshev polynomials of type II of degree $n_p$ in $\xi = \mathrm{arg}\,\nu_p$. The property (the second 
condition in \cref{Recursion}) of the coefficients in the $L$-series \cref{augDirichL} is the recursions relation \cref{Cheby2Recursion} 
satisfied by the Chebyshev polynomials. Finally, from the two equations above, we find that the coefficients are
\begin{equation}
a(n) = \sum_{d|n} \nu(d)\, \nu^*\!\left(\frac{n}{d}\right) = \prod_{\buildrel{p\: \text{such that}}\over{n=\prod p^{n_p}}} U_{n_p}
\left(\cos (\mathrm{arg}\, \nu_p)\right)       \label{CoeffAugDL}
\end{equation}
Before we proceed to discuss ${}_2\mathsf{L}(s,\nu)$ further, let us note that one may define a slightly more general product 
\begin{equation*}
{}_2\mathsf{L}(s,\nu_1,\nu_2^*) =  L(s,\nu_1) L(s,\nu_2^*) = \prod_p \left(1 - 2 e^{\frac{i}{2} \left(\alpha_1 - \alpha_2\right)}
\cos \frac{\alpha_1 + \alpha_2}{2}\,  p^{-s} +\chi(p) p^{-2s}\right)^{-1}
\end{equation*}
where $\alpha_1 = \mathrm{arg}\,\nu_1(p)$, $\alpha_2 = \mathrm{arg}\,\nu_2(p)$,  $\chi(p) = \nu_1(p) \nu_2^*(p) = e^{i\left(\alpha_1 - 
\alpha_2\right)}$ and the coefficients $a(p) = e^{\frac{i}{2} \left(\alpha_1 - \alpha_2\right)} U_1\left(\cos \frac{1}{2}(\alpha_1 + \alpha_2)\right)$ 
satisfy the multiplicative property $a(p) a(p^r) = a(p^{r+1}) + \chi(p) p^{k-1} a(p^{r-1})$ with $k=1$. This is a special case of the parametrisation
discussed at the end of \cref{ssec:IPtwo}.

The Dirichlet $L$-function in \cref{DirichL} is related to the $\vartheta$-series 
\begin{equation}
\vartheta(z,\nu) = \sum_{n\in\mathbb{Z}} \nu(n) n^\epsilon e^{i\pi n^2z/N}      \label{ThSeries}
\end{equation}
where $\nu$ is a primitive Dirichlet character modulo $N$ and $\epsilon = \frac{1}{2}(1-\nu(-1))$ takes the value $0$ or $1$ depending on 
whether $\nu$ is even or odd, respectively. The series above defines a modular form of weight $\frac{1}{2} + \epsilon$, level $4N^2$ and
character $\nu(d)\left(\frac{-1}{d}\right)^\epsilon$, where the Legendre symbol gives a phase. The Mellin transform of the above is the 
$L$-function \cref{DirichL}
\begin{equation}
L(s,\nu) =  \frac{\left({\pi}/{N}\right)^{\frac{s+\epsilon}{2}}}{2\Gamma\left(\frac{s+\epsilon}{2}\right)} \, \int_0^\infty \frac{dy}{y}
\, y^{\frac{s+\epsilon}{2}} \vartheta(iy, \nu)               \label{DLFnAsMellin}
\end{equation}
We may now use this in \cref{augDirichL} to write
\begin{align}
{}_2\mathsf{L}(s,\nu) &= \frac{({\pi}/{N})^{s+\epsilon}}{4\left(\Gamma\!\left(\frac{s+\epsilon}{2}\right)\right)^2} \, 
\int_0^\infty dy_1 \int_0^\infty dy_2\, (y_1 y_2)^{\frac{s+\epsilon}{2}-1} \vartheta(iy_1,\nu)\vartheta(iy_2,\nu^*)\nonumber\\
{} &= \frac{2({\pi}/{N})^{s+\epsilon}}{\left(\Gamma\!\left(\frac{s+\epsilon}{2}\right)\right)^2} \, \sum_{n_1=1}^\infty \nu(n_1) n_1^\epsilon 
\sum_{n_2=1}^\infty \nu^*(n_2) n_2^\epsilon \int_0^\infty \frac{dy}{y} y^{s+\epsilon}\, \int_0^\infty \frac{dy'}{y'}\, 
e^{- \frac{\pi n_1^2 yy'}{N} - \frac{\pi n_2^2 y}{Ny'}} \nonumber\\
{} &= \frac{4({\pi}/{N})^{s+\epsilon}}{\left(\Gamma\!\left(\frac{s+\epsilon}{2}\right)\right)^2} \, \sum_{n = 1}^\infty a(n) n^\epsilon 
\int_0^\infty \frac{dy}{y} y^{s + \epsilon}\, K_0\left(\frac{2\pi n}{N}\,y\right) \label{MellinAugDL}    
\end{align}
The integrals from the $L$-functions have been reorganised to be a convolution of the integrals of the two $\vartheta$-series. We have redefined 
the dummy variables of integrations and summations to $y^2 = y_1 y_2$, $y'^2 = y_1/y_2$, $n = n_1 n_2$, $d = n_1$ and used \cref{CoeffAugDL} 
for the coefficients. The integral in $y'$ can be evaluated from standard tables (e.g., from \cite{GradRyzh}), or identified as the integral 
representation of the modified Bessel function $K_0$. Finally, the $y$-integral in \cref{MellinAugDL}, being a Mellin transform of the modified 
Bessel function, is known analytically\cite{GradRyzh}. Performing the integral one of course gets back \cref{augDirichL} as expected.    

The expression above means that the function ${}_2\mathsf{L}(s,\nu)$ is the Mellin transform of the convolution of two $\vartheta$-series
\begin{align}
\left(\vartheta(\nu) \star \vartheta(\nu^*)\right)(iy) &= \frac{\left({\pi}/{N}\right)^{\frac{s+\epsilon}{2}}}{2\Gamma\left(\frac{s+\epsilon}{2}\right)} \, 
\int_0^\infty \frac{dy'}{y'} \vartheta(iyy',\nu)\, \vartheta\left(\frac{iy}{y'}, \nu^*\right)   \nonumber\\
&= \frac{\left({\pi}/{N}\right)^{\frac{s+\epsilon}{2}}}{2\Gamma\left(\frac{s+\epsilon}{2}\right)} \, \int_0^\infty \frac{dy'}{y'} \vartheta\left(\frac{iy}{y'},
\nu\right)\, \vartheta(iyy', \nu^*) \label{thetaConv} 
\end{align}
where the last equality follows from an exchange $n_1\leftrightarrow n_2$ along with a redefinition of the variable of integration 
$y'\rightarrow 1/y'$. The convolution is that of the Dirichlet character and its complex conjugate (denoted by a star above), as well as a Mellin 
convolution in the imaginary part of their arguments that is shown explicitly. We know that the $\vartheta$-series has modular property, in 
particular under $y \rightarrow 1/y$ (the $S$-transformation $z \rightarrow -1/z$ of the modular group restricted to the imaginary part) 
\begin{equation}
\vartheta\left(\frac{i}{y},\nu\right) = \frac{y^{\frac{1}{2}+\epsilon}}{i^\epsilon \sqrt{N}}\,\tau(\nu)\, \vartheta(iy, \nu^*) \label{thetaSTransf}
\end{equation}
where $\tau(\nu) = \displaystyle{\sum_{m=0}^{N-1}} \nu(m) e^{2\pi i m/N}$ is the Gauss sum. Hence we find that
\begin{equation}
\int_0^\infty \frac{dy'}{y'} \vartheta\left(\frac{i}{yy'},\nu\right)\star \vartheta\left(\frac{iy'}{y},\nu^*\right)    
= y^{1+2\epsilon}\,\int_0^\infty \frac{dy'}{y'}\, \vartheta\left(\frac{iy}{y'},\nu\right)\star \vartheta(iyy',\nu^*) \label{ConvThetaSTransf}
\end{equation}
where we have used $\tau(\nu)\tau(\nu^*) = \nu(-1)N$. Thus the Mellin inverse of ${}_2\mathsf{L}(s,\nu) $, namely the convolution 
\cref{thetaConv} of the $\vartheta$-series, is (quasi-)modular with weight $1+2\epsilon$ under the transformation $y \rightarrow 1/y$.      
  
It is interesting to note that this property, as well as the explicit expression on the RHS in \cref{MellinAugDL}, are reminiscent of the harmonic 
Maass waveform\cite{MaassZL} of weight $\lambda=0$ 
\begin{equation}
\mathsf{M}_{\lambda,N} (x+iy) = \sum_{n=1}^\infty a_n (n y)^\epsilon\, \sqrt{y}\, K_\lambda \left(\frac{2\pi n}{N}\, y\right) e^{2\pi i nx/N}  
\label{MaassForm}
\end{equation}
restricted to purely imaginary argument. A Maass waveform $\mathsf{M}_{\lambda,N} : \mathbb{H}\to \mathbb{C}$, is a non-holomorphic 
`modular function' on the upper half-plane that is a square-integrable eigenfunction of the hyperbolic $\Gamma(1)$-invariant Laplacian 
corresponding to the eigenvalue $\frac{1}{4} - \lambda^2$. Since $a_0=0$ in the above, it is actually a Maass cusp form. It may seem that there 
is a puzzle since the Maass waveform $\mathsf{M}_{0,N}$ has zero modular weight, while the $\vartheta$-series related to each of the 
$L$-functions \cref{DLFnAsMellin} are of weight $\frac{1}{2} + \epsilon$, hence their (convolution) product should be of weight $1+2\epsilon$. 
This is true, however, it is compensated by $y^{\frac{1}{2} + \epsilon}$, which is a non-holomorphic form of weight $-1-2\epsilon$. The factor of 
$\sqrt{y}$ was introduced by shifting the argument in the Mellin transform. It may also be noted that one loses complex analyticity by performing 
a Mellin transform in only the imaginary part of the argument of the $\vartheta$-series. 

We propose to identify the modular object related to the series ${}_2\mathsf{L}(s,\nu)$ to be a Maass form of the type above
\begin{equation}  
\sqrt{y}\, \mathbf{f}(x+iy,\nu) \rightarrow \mathsf{M}^{(\nu)}_{0,N} (x+iy) = \sum_{n=1}^\infty a(n) (n y)^\epsilon\,\sqrt{y}\, 
K_0\left(\frac{2\pi ny}{N}\right)\, e^{\frac{2\pi i nx}{N}}  \label{MaassCorres}
\end{equation}
where $a(n)$ is related to the Dirichlet character $\nu$ by \cref{CoeffAugDL}. We must caution that this identification is tentative since we have 
not been able to show the transformation property of \cref{thetaConv} under the full modular group. Even though this is only a special case of 
weight $\lambda=0$ and is moreover `reducible' by construction, this toy construction of (quasi-)Maass waveforms through a product of Dirichlet 
$L$-series could be useful to understand aspects of the proposed correspondence. The following analysis is in that spirit.

Henceforth we shall restrict to the case of even characters $\nu(1)=\nu(-1)$ for definiteness. The case of the odd characters may be studied in 
an analogous fashion. Let us now consider the $p$-artonic modular forms corresponding to \cref{augDirichL,GenCheby2,CoeffAugDL}
\begin{equation}
\mathbf{f}_{(p)}(\nu,x) = \sum_{n_p=0}^\infty a(p^{n_p}) \bupsi^{(p)}_{1-n_p,0,1}(x)
= \sum_{n_p=0}^\infty U_{n_p}\!\left(\cos(\mathrm{arg}\, \nu_p)\right) \bupsi^{(p)}_{1-n_p,0,1}(x) \label{augDirichLModF}  
\end{equation}
While these vectors in $\mathcal{H}_-^{(p)}$ are well-defined, the correspondence of the tensor product $\otimes_p \mathbf{f}_{(p)}(\nu,x_{(p)})$ 
to a modular object on the upper half-plane $\mathbb{H}$ proposed above would require further justification. This is because the Maass waveform 
\cref{MaassForm}, not being a meromorphic function, does not admit a $q$-series expansion. The latter form is what we had used for the 
correspondence \cref{qExpPFactorised} in \cref{sec:Parton}.
However, the dependence of both the holomorphic and non-holomorphic modular forms \cref{qExpand} and \cref{MaassCorres} on the variable 
$x$ is of the same form. Therefore, the association proposed should correctly be thought of as that between the Fourier coefficients in the Fourier 
series expansion (in $x$) of the corresponding modular objects. The $L$-functions related to the modular objects are also defined with the help of 
these Fourier coefficients. Nevertheless, it would be desirable to make a distinction between the holomorphic and non-holomorphic forms at the 
level of the proposed $p$-artons. 

Be that as it may, the inner product\footnote{For modular forms of weight $k=1$,
the two inner products \cref{pPeterssonIP} and \cref{pSimpleIP} are the same.} \cref{pSimpleIP} of two such functions 
$\mathbf{f}_{(p)}(\nu_{\mathbf{f}})$ and $\mathbf{g}_{(p)}(\nu_{\mathbf{g}})$ satisfies the orthogonality condition
\begin{align}
\big( \mathbf{f}_{(p)}(\nu_{\mathbf{f}}) \big| \mathbf{g}_{(p)}(\nu_{\mathbf{g}}) \big) &= \int_{\Qp^\times} d^\times x \; 
\mathbf{f}_{(p)}^*(\nu_{\mathbf{f}},x)\, \mathbf{g}_{(p)}(\nu_{\mathbf{g}}, x) \nonumber\\
&= \sum_{n_p=0}^\infty U_{n_p}\!\left(\cos(\mathrm{arg}\, \nu_{\mathbf{f},p}^*)\right)\, U_{n_p}\!\left(\cos(\mathrm{arg}\, \nu_{\mathbf{g},p})\right) 
\,=\, \delta_{\nu_{\mathbf{f}},\nu_{\mathbf{g}}} \label{OrthoAugDirichLModF}
\end{align}
as a consequence of the completeness of the Chebyshev polynomials. Hence the inner product of their tensor products 
$\mathbf{f}(\nu_{\mathbf{f}})$ and $\mathbf{g}(\nu_{\mathbf{g}})$
\begin{align}
\big( \mathbf{f}(\nu_{\mathbf{f}}) \big| \mathbf{g}(\nu_{\mathbf{g}}) \big) &= \prod_p \big( \mathbf{f}_{(p)}(\nu_{\mathbf{f}}) \big| 
\mathbf{g}_{(p)}(\nu_{\mathbf{g}}) \big) \nonumber\\
&= \prod_p \sum_{n_p=0}^\infty a^*_{\mathbf{f}}(p^{n_p})\, a_{\mathbf{g}}(p^{n_p}) \nonumber\\
&= \sum_{n=1}^\infty a^*_{\mathbf{f}}(n)\, a_{\mathbf{g}}(n) \label{fgIP1}
\end{align}
will also be orthogonal. Let us consider the inner product \cref{PeterssonIP} of two Maass forms $\mathsf{M}_{\mathbf{f}} = \sqrt{y}\,
\mathbf{f}(x+iy,\nu_{\mathbf{f}})$ and $\mathsf{M}_{\mathbf{g}} = \sqrt{y}\,\mathbf{g}(x+iy,\nu_{\mathbf{g}})$
\begin{equation*}
\big\langle \mathbf{f}(\nu_{\mathbf{f}}) \big| \mathbf{g}(\nu_{\mathbf{g}}) \big\rangle = \int_{{-\frac{N}{2}<|x|\le\frac{N}{2}}\atop{|x+iy|>0}} 
\frac{dx\,dy}{y^2} \left(\sqrt{y}\,\mathbf{f}(x+iy,\nu_{\mathbf{f}})\right)^*\,\sqrt{y}\,\mathbf{g}(x+iy,\nu_{\mathbf{g}}) 
\end{equation*}
Unfortunately, we are not able to perform the integrals with the restrictions imposed by the fundamental domain. However, we can compute the
above in the rectangular region $\{-\frac{N}{2}<|x|\le\frac{N}{2},\, y > 0\}$, which contains an {\em infinite} number of copies of the fundamental
domain $\mathcal{F}$. By the modular properties of the integrand and the measure, each of these copies are identical, therefore, we expect
to get the original integral multiplied by an infinite factor. That is exactly what we find. The $x$-integral in
\begin{align}
\big\langle \mathbf{f}(\nu_{\mathbf{f}}) \big| \mathbf{g}(\nu_{\mathbf{g}}) \big\rangle &= \sum_{m,n=1}^\infty a^*_{\mathbf{f}}(m)
a_{\mathbf{g}}(n) \int_{-{N}/{2}}^{{N}/{2}} dx\, e^{2\pi i(n-m)x/N} \int_0^\infty \frac{dy}{y} K_0\left(\frac{2\pi my}{N}\right)\,
K_0\left(\frac{2\pi ny}{N}\right)\nonumber 
\end{align}
vanishes unless $m=n$, giving $\delta_{mn}$, implying that two different Maass forms are orthogonal. For the $y$-integral we use the 
expression \cite{GradRyzh}
\begin{align*}
\int_0^\infty dy\, y^{-\mu} K_\alpha(ay)\, K_\beta(by) &= \frac{a^{\mu-1-\beta}b^\beta}{2^{\mu+2}\Gamma(1-\mu)}\,
\Gamma\left(\frac{1-\mu+\alpha+\beta}{2}\right) \Gamma\left(\frac{1-\mu-\alpha+\beta}{2}\right)\\
&\quad\times\Gamma\left(\frac{1-\mu+\alpha-\beta}{2}\right) \Gamma\left(\frac{1-\mu-\alpha-\beta}{2}\right)\\
&\quad\times {}_2F_1\left(\frac{1-\mu+\alpha+\beta}{2}, \frac{1-\mu-\alpha+\beta}{2},1-\mu;1-\frac{b^2}{a^2}\right)\\
&\quad\qquad\text{for }\mathrm{Re}\,(a+b) > 0,\: \mathrm{Re}\,\mu < 1 - |\mathrm{Re}\,\alpha| - |\mathrm{Re}\,\beta| 
\end{align*}
The integral we need is exactly at the border of the condition on $\mu$, and indeed the arguments of all the $\Gamma$-functions are 
zero, giving a divergent factor. Thus, upto an infinite normalisation factor arising out of the divergent $\Gamma$-functions, the modular 
objects related to the product $L$-functions ${}_2\mathsf{L}(s,\nu)$ are orthogonal non-analytic Maass waveforms of modular weight zero. 
This corresponds to the orthogonality \cref{fgIP1} obtained from the wavelet expansion.  

The expression \cref{fgIP1} above can be related to the inner product of the corresponding product $L$-functions \cref{augDirichL} by using 
the identity
\begin{equation*}
\lim_{T\to\infty} \frac{1}{2T} \int_{-T}^{T} dt\, \frac{e^{i t\left(\ln m - \ln n\right)}}{\sqrt{mn}} = \frac{\delta\left(\ln n - \ln m\right)}{\sqrt{mn}} 
= \delta_{mn}
\end{equation*}
It should be mentioned that while the identity above for $\delta_{mn}$ works for $m^\alpha n^{1-\alpha}$ for {\em any} $0<\alpha<1$ in the 
denominator, only for $\alpha = \frac{1}{2}$ does it lead to the properties required of an inner product in the following. Thus
\begin{align}
\big\langle \mathbf{f}(\nu_{\mathbf{f}}) \big| \mathbf{g}(\nu_{\mathbf{g}}) \big\rangle  
&= \sum_{n=1}^\infty \sum_{m=1}^\infty a^*_{\mathbf{f}}(n)\, a_{\mathbf{g}}(m) \nonumber\\
&= \lim_{T\to\infty} \frac{1}{2T} \int_{-T}^{T} dt\, \sum_{n=1}^\infty \frac{a^*_{\mathbf{f}}(n)}{m^{\frac{1}{2}-it}} \, 
\sum_{n=1}^\infty \frac{a_{\mathbf{g}}(n)}{n^{\frac{1}{2}+it}}\nonumber\\
&= \lim_{T\to\infty} \frac{1}{2T} \int_{-T}^{T} dt\; {}_2\mathsf{L}_{\mathbf{f}}^*\left(\frac{1}{2} - it, \nu_{\mathbf{f}} \right) \, 
{}_2\mathsf{L}_{\mathbf{g}}\left(\frac{1}{2} + it, \nu_{\mathbf{g}}\right)  \label{fLgLIP}
\end{align}
where in the last line we have used the analytical continuation of the power series. The final expression allows for an interpretation as an inner 
product $\big\langle {}_2\mathsf{L}_{\mathbf{f}} \big| {}_2\mathsf{L}_{\mathbf{g}} \big\rangle$ of the $L$-functions related to the Maass forms. 
This is a Parseval type identity, like that in \cref{ParsevalType,modParsevalType}. Hence, the family of product $L$-functions, defined in 
\cref{GenCheby2}, seems to form an orthogonal set of functions (for a fixed prime value of $N$). Although the $L$-functions are known explicitly, 
it is not straightforward to verify their orthogonality. It may be seen (by plotting the functions in {\tt Mathematica}) that the imaginary part of 
the integrand is odd, therefore, its integral vanishes for $\mathbf{f}\ne\mathbf{g}$, however, its real part oscillates with increasing amplitude 
making it difficult to perform the integral numerically. 

\section{Endnote}\label{sec:Summary}
In this paper we have used the orthonormal bases provided by the wavelets \cite{Kozyrev:2001} in $L^2(\Qp)$ to associate complex valued 
functions, one function for each prime number $p$ to a cusp form of a congruence subgroup of the modular group SL(2,$\mathbb{Z})$. We 
have studied the {\em local} functions (which may appropriately be called $p$-artons, borrowing a term from the physics of strong interaction) 
obtained through this decomposition, their Mellin transforms and the related $L$-functions. We leave the study of the action of the group 
GL(2,$\Qp$) on these $p$-artons to the future. By taking a product of two Dirichlet $L$-functions, associated to a Dirichlet character and its
complex conjugate, we have also defined functions that are similar to these modular $L$-functions. Rather surprisingly these turn out to 
behave like non-analytic Maass waveforms of weight zero. Although being reducible by construction, these may not be of intrinsic interest 
mathematically, however, being much simpler to work with, they may be useful toy objects to further the correspondence. 

Let us close with the following remark. The one to one relation between a cusp form (a complex function on the upper half of the complex
plane) and a vector in (a subspace of) $\otimes_p L^2(\Qp)$ is reminiscent of the holographic correspondence, which has dominated the 
landscape of research in theoretical physics in the recent decades. The conformal boundary of the upper half plane 
($\mathbb{H} = \mathrm{SL}(2,\mathbb{R})/\mathrm{U}(1)$) is the real line $\mathbb{R}$. The latter is closely related to $\Qp$, which can 
be thought of as the conformal boundary of the Bruhat-Tits tree, which in turn is the coset of GL(2,$\Qp$) by its maximal compact subgroup
GL(2,$\Zp$). The association between a modular form $f:\mathbb{H}\to\mathbb{C}$ and the products of functions $\otimes_p\mathfrak{f}_{(p)}$,  
where $\mathfrak{f}_{(p)} : \Qp\to\mathbb{C}$ for a prime $p$ is a complex valued function on $\Qp$, is suggestive of a bulk-boundary 
correspondence in holography. It is so, since the data related to a function in the bulk of the upper half of the complex plane is in the `boundary'
$\otimes_p\Qp$, which in turn is related to $\mathbb{R}=\partial\mathbb{H}$. In this sense, the proposed relation may even be called a {\em 
holographic $p$-arton model} of modular forms. It would be interesting to explore if there are conformal field theories (and their bulk duals) 
related to the $p$-atronic objects as well as the non-analytic Maass waveforms.   

\bigskip

\noindent{\bf Acknowledgments:} DG is supported in part by the Mathematical Research Impact Centric Support (MATRICS) grant no.\
MTR/2020/000481 of the Science \&\ Engineering Research Board, Department of Science \&\ Technology, Government of India. 
One of us (DG) presented a preliminary version of some of these results (as well as those in \cite{Dutta:2020qed}) in the \emph{National 
String Meeting 2019} held at IISER Bhopal, India during Dec 22--27, 2019. We would like to thank Chandan Singh Dalawat, and Vijay Patankar 
for useful discussions. We are particularly grateful to Krishnan Rajkumar for many patient explanations and valuable comments on the manuscript. 

 
\appendix
\section{Wavelets on $\Qp^\times$}\label{app:AppA}
Let us form modify the Kozyrev wavelets \cref{waveket} to define
\begin{equation}
\bupsi^{(p)}_{n,m,j}(x) = |x|_p^{\frac{1}{2}}\, \psi^{(p)}_{n,m,j}(x)  \label{pWaveletMult} 
\end{equation}  
which are naturally defined on $\Qp^\times$ since they are orthonormal with respect to the scale invariant multiplicative Haar measure 
$d^\times x$:
\begin{equation*}
\int_{\Qp^\times} \frac{dx}{|x|_p} \bupsi^{(p)}_{n,0,1}(x) \bupsi^{(p)}_{n',0,1}(x) =
\int_{\Qp} dx\, \psi^{(p)}_{n,0,1}(x) \psi^{(p)}_{n',0,1}(x) = \delta_{nn'}  \label{modOrthonorm}
\end{equation*}
Notice that the wavelets above, being different from the Kozyrev wavelets \cref{KozyWvlet,waveket} by a coordinate dependent factor, is 
not the constant $p^{-n/2}$ for $|x|_p < p^n$, although it is still a locally constant function on $\Qp$. The raising and lowering operators
\cref{raise:lower,annihilate} act as before: $\mathbf{a}^{(p)}_\pm\bupsi^{(p)}_{n,0,1}(x)  = \bupsi^{(p)}_{n\pm 1,0,1}(x)$
and $\mathbf{a}^{(p)}_+ \bupsi^{(p)}_{1,0,1}(x) = 0$, however, we choose to use a different notation to emphasise the fact that they
act on a different space of functions.  The Mellin transform \cref{pMellin} of these modified wavelets 
\begin{equation}
\mathcal{M}_{(p,\omega)}[\bupsi^{(p)}_{n,0,1}](s) = \mathbf{c}_p(\ell,s)\, p^{ns} = - \left(\frac{1}{p\left(1 - p^{-s-\frac{1}{2}}\right)} - 
\frac{1}{p^{s+\frac{1}{2}} - 1}\,\delta_{\ell,0} - \delta_{\ell,p-1}\right) p^{ns}\label{MellinModWvlt}
\end{equation}
likewise differ somewhat from \cref{KozyMellin}. Let us note that $\sum_{\ell} |\mathbf{c}_p(\ell,s)|^2 = 1 + 
\frac{1 - |p^s|^2}{\big| p^{s+\frac{1}{2}} - 1\big|^2}$, thus, if the argument $s$ is purely imaginary, the sum is 1.

\bigskip

\bibliographystyle{hieeetr}
\bibliography{PseuD}{}

\end{document}